\documentclass[11pt]{article} %draft

\usepackage{amsmath, amsfonts, amssymb}
\usepackage{mathrsfs}
\usepackage{theorem}
\usepackage{bm}
\usepackage{xcolor}
\usepackage{anyfontsize}
\usepackage{tabularx}
\usepackage{colortbl}
\usepackage{graphicx}
\usepackage{hyperref}%
\usepackage{upgreek}
\usepackage{verbatim}
\usepackage{array}
\usepackage{graphicx}
\usepackage{tikz}

\newtheorem{mythm}{Theorem}[section]
\newtheorem{myprop}[mythm]{Proposition}
\newtheorem{mylem}[mythm]{Lemma}
\newtheorem{mycor}[mythm]{Corollary}
\newtheorem{mydefn}[mythm]{Definition}

{\newtheorem{myrem}[mythm]{Remark}}
{%\theorembodyfont{\rmfamily}
\newtheorem{myexam}[mythm]{Example}}%
\numberwithin{equation}{section}

\allowdisplaybreaks

\begin{document}

\title{ Non-(strong, geometrically) ergodicity criteria for discrete time Markov chains on general state spaces\footnote{$*$ Corresponding author: Yu Chen; E-mail: ychen0504@mail.bnu.edu.cn} }
\author{ Ling-Di Wang$^a$\thanks{$a$. School of Mathematics and Statistics, Henan University, Kaifeng 475001, China;}\quad Yu Chen$^{b*}$\thanks{$b$. School of Mathematics Sciences, Beijing Normal University, Laboratory of Mathematics and Complex Systems, Ministry of Education, Beijing 100875, China.} \quad Yu-Hui Zhang$^b$
}

%\date{}
\maketitle

\begin{abstract}For discrete-time Markov chains on general state spaces, we establish criteria for non-ergodicity and non-strong ergodicity, and derive sufficient conditions for non-geometric ergodicity via the theory of minimal nonnegative solutions. Our criteria are formulated based on the existence of solutions to inequalities involving the chain’s one-step transition kernel. Meanwhile, these practical criteria are applied to a type of examples, which can effectively characterize the non-ergodicity and non-strong ergodicity of a specific class of single birth (death) processes.
\end{abstract}



\section{Introduction}Ergodicity is one of the three fundamental and classical problems (uniqueness, recurrence, and ergodicity) for Markov chains (\cite{2004-Chen}). Ordinary ergodicity, strong ergodicity, and so on, are usually used to characterize the stability of Markov chains, which play important roles in both theoretical and applied research (\cite{2002-J-R, 1995-L, 1992-M-T, 1993-M-T, 1990-T}). The topic of criteria for various types of ergodicity has been studied by many authors over the past decades (see \cite{2004-Chen, 2003-Mao, 2014-M-S, 2009-M-T}, etc.). In particular, Lyapunov (drift) criteria have long been used to provide necessary and sufficient conditions for ergodicity or strong ergodicity. In principle, these results also imply conditions for non-ergodicity and non-strong ergodicity by stating the non-existence of such Lyapunov functions, which is usually not practical. For this reason, criteria for non-ergodicity and non-strong ergodicity of Markov chains involving test functions—known as the "inverse problem"—have also been specially studied. For the inverse problem of Markov processes, we refer readers to \cite{2004-C-K, 1983-S-T, 1985-S} for denumerable state spaces and continuous-time settings, to \cite{2022-Wei} for denumerable state spaces and discrete-time settings, and to \cite{2021-M-W} for general state spaces and continuous-time settings. However, to our knowledge, there are few results in the existing literature on non-ergodicity or non-strong ergodicity for discrete-time Markov chains on general state spaces. This gap is addressed in the present paper, where we establish criteria and sufficient Lyapunov (drift) conditions for their non-ergodicity, non-strong ergodicity and non-geometric ergodicity.

Consider a homogeneous Markov chain on a general state space $\mathscr{X}$, which is a locally compact separable metric topological space endowed with a countably generated $\sigma$-field $\mathscr{B}(\mathscr{X})$. Denote by
\[
P^n(x,A) = \mathbb{P}(X_n \in A \mid X_0 = x), \quad n \ge 0,\ x \in \mathscr{X},\ A \in \mathscr{B}(\mathscr{X}),
\]
with $P(x, A) = P^1(x, A)$.
The Markov chain $\{X_n\}$ is called $\psi$-irreducible if there exists a $\sigma$-finite measure $\psi$ on $(\mathscr{X}, \mathscr{B}(\mathscr{X}))$ such that $\sum_{n=1}^{\infty} P^n(x, A) > 0$ for all $x \in \mathscr{X}$ whenever $\psi(A) > 0$.

\textbf{Throughout the paper, we assume the chain $\{X_n\}$ is $\psi$-irreducible, where $\psi$ is a maximal irreducibility measure on $(\mathscr{X}, \mathscr{B}(\mathscr{X}))$.}

The concepts of "aperiodicity", "Feller", and "$\psi$-irreducible" coincide with those in \cite{2009-M-T} (see Section 5.4.2 and Chapter 6 therein), and we omit their definitions here for simplicity. Let $Pf(x) = \int_{\mathscr{X}} f(y) P(x, \mathrm{d}y)$ for some $\mathscr{B}(\mathscr{X})$-measurable function $f$, and define
\[
\mathscr{B}^+(\mathscr{X}) = \left\{ A \in \mathscr{B}(\mathscr{X}) : \psi(A) > 0 \right\}.
\]
If there exists a probability distribution $a = \{a_n\}$ on $\mathbb{Z}_+$ such that
\[
\sum_{n=0}^{\infty} a_n P^n(x, B) \ge \nu_a(B), \quad B \in \mathscr{B}(\mathscr{X})
\]
for all $x \in A$, where $\nu_a$ is a non-trivial measure on $\mathscr{B}(\mathscr{X})$, then $A \in \mathscr{B}(\mathscr{X})$ is called a $\nu_a$-petite set (see \cite[Section 5]{2009-M-T}).
Petite sets are not rare: compact sets are petite under some mild conditions \cite[Proposition 6.2.8]{2009-M-T}. For each $A \in \mathscr{B}^+(\mathscr{X})$, there exists a petite set $B \subseteq A$ with $B \in \mathscr{B}^+(\mathscr{X})$ \cite[Theorem 5.2.2]{2009-M-T}.

The Markov chain $\{X_n\}$ is called 
\begin{itemize}
\item[1)] \textbf{recurrent} if 
$\sum_{n=1}^{\infty} P^n(x,A) = \infty$ for all $x \in \mathscr{X}$ and $A \in \mathscr{B}^+(\mathscr{X})$; otherwise, it is called \textbf{transient};

\item[2)] \textbf{Harris recurrent} if $\mathbb{P}_x\left( \sum_{n=1}^{\infty} \mathbf{1}_{\{X_n \in A\}} = \infty \right) = 1$ for all $x \in \mathscr{X}$ and $A \in \mathscr{B}^+(\mathscr{X})$;

\item[3)] \textbf{positive} if it admits an invariant probability measure. Equivalently, for any subinvariant measure $\mu$, and for one (hence every) set $A$ with $\mu(A) > 0$,
\[
\int_A \mathbb{E}_y \tau_A^+ \, \mu(\mathrm{d}y) < \infty;
\]

\item[4)] \textbf{ergodic} if it is recurrent and there exists some petite set $C \in \mathscr{B}^+(\mathscr{X})$ such that $\sup_{x \in C} \mathbb{E}_x \tau_C^+ < \infty$. Equivalently, it is positive and Harris recurrent;

\item[5)] \textbf{geometrically ergodic} if there exists some petite set $C \in \mathscr{B}(\mathscr{X})$ and $r > 1$ such that 
$\sup_{x \in C} \mathbb{E}_x r^{\tau_C^+} < \infty$;

\item[6)] \textbf{strongly (or uniformly) ergodic} if there exists a petite set $C \in \mathscr{B}(\mathscr{X})$ with 
$\sup_{x \in \mathscr{X}} \mathbb{E}_x \tau_C^+ < \infty$. In this case, for every set $A \in \mathscr{B}^+(\mathscr{X})$, $\sup_{x \in \mathscr{X}} \mathbb{E}_x \tau_A^+ < \infty$.
\end{itemize}

For more details on these notions and their equivalent descriptions, refer to \cite{2009-M-T, 1990-T}. In general, for aperiodic Harris recurrent Markov chains, "positive Harris recurrent chains" are referred to as "ergodic chains" (\cite[Theorem 13.0.1]{2009-M-T}). A Harris recurrent chain is necessarily recurrent. Since we focus on ergodic properties, we assume that all Markov chains considered below are Harris recurrent unless specified otherwise. In this paper, "Harris recurrent" is abbreviated as "recurrent".

Now, we present our results.

\textbf{Assumption 1}
The Markov chain $\{X_n\}$ is $\psi$-irreducible, aperiodic, Feller, and the support of $\psi$ (i.e., $\mathrm{supp}\,\psi$) has nonempty interior.

\begin{mythm}\label{Non-1-ergo}(Non-ergodicity) Assume that the Markov chain $\{X_n\}$ satisfies Assumption 1. Then $\{X_n\}$ is not ergodic if there exist a set $A \in \mathscr{B}^+(\mathscr{X})$ and a sequence of functions $\{V^{(n)}(x)\}_{n \ge 1}$ with $V^{(n)}(x): \mathscr{X} \to \mathbb{R}$ that satisfy the following conditions:
\begin{itemize}
\item[1)] $\sup_{x \in A^c} V^{(n)}(x) < \infty$ for all $n \ge 1$;
\item[2)] $\int_{A^c} V^{(n)}(y) P(x, \mathrm{d}y) \ge V^{(n)}(x) - 1$ for all $x \in \mathscr{X}$ and $n \ge 1$;
\item[3)] $\int_{A} \sup_{n \ge 1} V^{(n)}(y) \psi(\mathrm{d}y) = \infty$.
\end{itemize}
The converse also holds when the Markov chain $\{X_n\}$ is recurrent.
\end{mythm}

\begin{mycor}\label{Non-1-ergo-cor} Assume that the Markov chain $\{X_n\}$ satisfies Assumption 1 and condition 1) of Theorem \ref{Non-1-ergo}. The chain $\{X_n\}$ is not ergodic if inequality 2) of Theorem \ref{Non-1-ergo} holds for all $x \in \mathscr{X}$ (or for all $x \in A^c$) and 
\[
\left\{ x \in B : \sup_{n \ge 1} V^{(n)}(x) = \infty \right\} \in \mathscr{B}^+(\mathscr{X})
\]
for some set $B \in \mathscr{B}^+(\mathscr{X})$ (or for $B \in \mathscr{B}^+(\mathscr{X})$ with $B \subset A^c$, respectively).
\end{mycor}

\begin{myrem} Let $\mu$ denote the invariant measure when the chain is recurrent. Then $\mu$ and $\psi$ are equivalent. Since $P(x, \cdot) \ll \mu$ for $\mu$-almost all $x$, we have $P(x, \cdot) \ll \psi$ for $\mu$-almost all $x$. Thus, the sufficiency condition in Theorem \ref{Non-1-ergo} remains valid if condition 3) is replaced by
\[
\int_{A} \sup_{n \ge 1} V^{(n)}(y) P(x, \mathrm{d}y) = \infty \quad \text{for } \psi\text{-almost all } x.
\]
\end{myrem}\begin{mythm}\label{Non-stro}(Non-strong ergodicity) Assume that the Markov chain $\{X_n\}$ satisfies Assumption 1. Then $\{X_n\}$ is non-strongly ergodic if and only if there exist a set $A \in \mathscr{B}^+(\mathscr{X})$ and a sequence of functions $\{V^{(n)}(x)\}_{n \ge 1}$ with $V^{(n)}(x): \mathscr{X} \to \mathbb{R}$ that satisfy the following conditions:
\begin{itemize}
\item[1)] $\sup_{x \in A^c} V^{(n)}(x) < \infty$ for all $n \ge 1$; and $V^{(n)}(x) = 0$ for all $x \in A$ and $n \ge 1$;
\item[2)] $PV^{(n)}(x) \ge V^{(n)}(x) - 1$ for all $x \in A^c$ and $n \ge 1$;
\item[3)] $\sup_{\substack{x \in A^c \\ n \ge 1}} V^{(n)}(x) = \infty$.
\end{itemize}
\end{mythm}

\begin{myrem} If we replace Assumption 1 with "the Markov chain $\{X_n\}$ is $\psi$-irreducible and aperiodic" in Theorem \ref{Non-1-ergo} and Theorem \ref{Non-stro}, their sufficiency conditions still hold, respectively.
\end{myrem}

\begin{myprop}\label{Dykin} Assume that $\{E_n, n \ge 1\} \subset \mathscr{B}(\mathscr{X})$ is an increasing sequence of bounded sets with $\mathscr{X} = \cup_{n=1}^{\infty} E_n$. Suppose there exists a continuous function $V: \mathscr{X} \to \mathbb{R}$ and a set $A \in \mathscr{B}^+(\mathscr{X})$ such that
\begin{itemize}
\item[1)] $\sup_{x \in A} V(x) < \infty$ and $\sup_{x \in A^c} V(x) = \infty$;
\item[2)] $V(x)$ is locally bounded, and $PV(x) \ge V(x) - 1$ for all $x \in A^c$;
\item[3)] $W(x): \mathscr{X} \to [0, \infty)$ is continuous, $PW(x) \le W(x) + d\mathbf{1}_{A}(x)$ for some constant $d$, and 
\[
\lim_{m \to \infty} \sup_{x \in E_m^c} \frac{V(x)}{W(x)} = 0.
\]
\end{itemize}
Then $\{X_n\}_{n \ge 0}$ is non-strongly ergodic.
\end{myprop}

\begin{mythm}\label{Non-geo}(Non-geometric ergodicity) Assume that the Markov chain $\{X_n\}$ is $\psi$-irreducible and aperiodic. Then $\{X_n\}$ is non-geometrically ergodic if there exist a set $A \in \mathscr{B}^+(\mathscr{X})$, a sequence of numbers $\{r_n\}_{n=1}^{\infty}$ with $r_n > 1$, and a sequence of functions $\{V^{(n)}(x)\}_{n \ge 1}$ with $V^{(n)}(x): \mathscr{X} \to \mathbb{R}$ that satisfy the following conditions:
\begin{itemize}
\item[1)] $\lim_{n \to \infty} r_n = 1$;
\item[2)] $V^{(n)}(x)$ is compactly supported for all $n \ge 1$, bounded on any compact set, and solves the inequality
\[
V^{(n)}(x) \le r_n \int_{A^c} V^{(n)}(y) P(x, \mathrm{d}y) + 1 \quad \text{for all } x \in \mathscr{X} \text{ and } n \ge 1;
\]
\item[3)] $\{ x : \sup_{n \ge 1} V^{(n)}(x) = \infty \} \in \mathscr{B}^+(\mathscr{X})$.
\end{itemize}
\end{mythm}
Note that in the proof of the above criteria, we reveal an interesting ``symmetry''. To be brief, various types of moments of hitting times, as minimal solutions to some equations, are also nicely bounded below by solutions to some corresponding ``symmetric'' inequalities and therefore exhibit some ``maximal'' property. Although most of our criteria need a sequence of testing functions, which makes it seemingly ``involved'', we can actually do batch production for testing functions. For example, one may consult the following proposition, and its proof of the generalized version in Section 3.1 (see Propositions \ref{ex1-transi} $\sim$ \ref{ex1-str-erg}), which is rather easy using the assertions obtained in the paper (Theorem \ref{Non-1-ergo} and Theorem \ref{Non-stro}, Proposition \ref{transient-1}, etc.) combined with the existing sufficient conditions for ergodicity and related properties.  in \cite{2009-M-T}.

\begin{myprop}\label{1-Ex-1}
Let $1>\gamma(x)>0$ be a continuous function on $\mathbb{R}^+$, and let $\beta(x)$ be a density function with support set $(0, a)$,  $0<a\le \infty$. $(X_n)_{n\ge0}$ is a Markov chain on $\mathbb{R}^+$ with one-step transition kernel $P(x, {\rm d}y)$ given by
\begin{equation*}
P(x,{\rm d}y)=\begin{cases}
\beta(y){\rm d}y,&x=0\\
\gamma(x)\delta_{\{0\}}({\rm d}y)+(1-\gamma(x))\delta_{\{x+1\}}({\rm d}y), &x>0.
\end{cases}
\end{equation*}
Then
\begin{itemize}
\item[1)] $\{X_n\}$ is recurrent if there exists $x_1>0$ such that $$\prod_{k=[x_1]}^{\infty}\sup_{k\le x<k+1}(1-\gamma(x))=0.$$ $\{X_n\}$ is transient if there exists $x_1>0$ such that 
$$\prod_{k=[x_1]}^{\infty}\inf_{k\le x<k+1}(1-\gamma(x))>0.$$
\item[(2)] $\{X_n\}$ is non-ergodic if there exists $x_1>0$ such that $$\int_{x_1}^{\infty}\frac{\beta(x)}{\sup_{y\ge x}\gamma(y)}{\rm d}x=\infty.$$ $(X_n)$ is ergodic if there exists $x_1>0$ such that 
$$\sup_{x>x_1}\left(\frac{1}{\gamma(x+1)}-\frac{1}{\gamma(x)}\right)<1,
\qquad \int_{x_1}^{\infty}\frac{\beta(y)}{\gamma(y)}{\rm d}y<\infty.$$
\item[(3)] $\{X_n\}$ is non-geometrically ergodic and non-strongly ergodic if $\overline{\lim}_{x\to\infty}\gamma(x)=0$.
\item[(4)] $\{X_n\}$ is strongly ergodic if $\underline{\lim}_{x\to\infty}\gamma(x)>0$.
\end{itemize}
\end{myprop}\begin{myrem}
Note that although explicit criteria for ergodicity and strong ergodicity are already given for such examples (e.g., the continuous-time single birth process with immigration in a countable state space), the above proposition would be difficult to derive without our criteria. After all, we only have explicit sufficient conditions for the exponential ergodicity of single birth processes with immigration in \cite{2010-Zhang-Zhao}. Using the above proposition, we present the following interesting assertions. 
\begin{itemize}
\item Assume that $a>1$. Let $\gamma(x)=1-(\sin^2x)/a$, $x\ge0$. Then $\{X_n\}$ is strongly ergodic since $\underline{\lim}_{x\to\infty}\gamma(x)=1-1/a>0$.
\item Assume that $r>0$, and let $\gamma(x)=x^{-r}$ for $x\ge1$ be a continuous function. 
\begin{itemize} 
\item[1)] $\{X_n\}$ is recurrent when $0<r\le 1$, and transient when $r>1$. 
\item[2)] Assume that $0 < \delta \leq r \le1$,
\[
\beta(x) =\frac{\delta}{\delta + 1}\left[{\bf 1}_{[0,1]}(x)+\frac{1}{x^{\delta+1}}{\bf 1}_{(1,\infty)}(x)\right]
\]
Then $\{X_n\}$ is null recurrent (recurrent but not ergodic).
\item[3)] Assume that $0 < r < 1$, and $\beta(x)$ is a density function of a random variable $\xi$ having small moments (i.e., $\mathbb{E}\xi^r<\infty$, for example, the Uniform distribution). Then $\{X_n\}$ is ergodic. However, it is non-geometrically ergodic and non-strongly ergodic.
\end{itemize}
\end{itemize}
\end{myrem}

\section{Proof of criteria for inverse problems}
\subsection{Recall minimal solution theory}
We first recall some useful results from the minimal solution theory in \cite{2004-Chen}.

Define $\mathscr{H}$ as a set of mappings from $\mathscr{X}$ to $[0, \infty]$: $\mathscr{H}$ contains the constant function $1$ and is closed under nonnegative combinations and monotone increasing limits, where the order relation ``$\ge$'' in $\mathscr{H}$ is defined point-wise. Then $\mathscr{H}$ is a convex cone. A mapping $H: \mathscr{H} \to \mathscr{H}$ is called a cone mapping if $H0 = 0$ and, for any nonnegative constants $c_1, c_2$ and mappings $f_1, f_2 \in \mathscr{H}$, 
\[
H(c_1f_1 + c_2f_2) = c_1H(f_1) + c_2H(f_2).
\]
Denote by $\mathscr{A}$ the set of all such cone mappings that also satisfy the following hypothesis: for $f_n \in \mathscr{H}$,
\[
f_n \uparrow f \implies Hf_n \uparrow Hf.
\]

\begin{mydefn} (\cite[Definition 2.1]{2004-Chen}) Given $H \in \mathscr{A}$ and $g \in \mathscr{H}$, we say $f^*$ is the minimal nonnegative solution (abbreviated as minimal solution) to the equation
\begin{equation}\label{mini-solution}
f(x) = (Hf)(x) + g(x), \quad x \in \mathscr{X},
\end{equation}
if $f^*$ satisfies \eqref{mini-solution} and for any solution $\widetilde{f} \in \mathscr{H}$ of \eqref{mini-solution},
\[
\widetilde{f}(x) \ge f^*(x), \quad x \in \mathscr{X}.
\]
This last property is called the minimality property of $f^*$.
\end{mydefn}

\begin{mylem}\label{mini-exsits} (\cite[Theorem 2.2]{2004-Chen}) The minimal solution to Equation \eqref{mini-solution} always exists uniquely.\end{mylem}

By Lemma \ref{mini-exsits}, for each $H \in \mathscr{A}$, we may define a map $m_H$ from $\mathscr{H}$ into itself by $m_H(g) = f^*$. We note that for $H_1, H_2 \in \mathscr{A}$, $H_1 \ge H_2$ if and only if $H_1f \ge H_2f$ for all $f \in \mathscr{H}$.

\begin{mylem}\label{comparison} (\cite[Theorem 2.6]{2004-Chen}) Let $f^*$ be the minimal solution to \eqref{mini-solution}. Then for any solution $\tilde{f}$ to 
\begin{equation}\label{mini-solution-control}
\tilde{f}(x) \ge (\tilde{H} \tilde{f})(x) + \tilde{g}(x), \quad x \in \mathscr{X},
\end{equation}
where $\tilde{H} \in \mathscr{A}$, $\tilde{g} \in \mathscr{H}$ with $\tilde{H} \ge H$ and $\tilde{g} \ge g$, we have $\tilde{f}(x) \ge f^*(x)$ for all $x \in \mathscr{X}$.
\end{mylem}

\begin{mylem}\label{Monot} (\cite[Theorem 2.7]{2004-Chen}) $m_H$ is a cone mapping. Let $\{g_n\}_{n=1}^{\infty} \subset \mathscr{H}$ and $\{H_n\} \subset \mathscr{A}$ with $H_n \uparrow H$ and $g_n \uparrow g$. Then $g \in \mathscr{H}$, $H \in \mathscr{A}$, and $m_{H_n}(g_n) \uparrow m_H(g)$.
\end{mylem}

\begin{mylem}\label{Second-ite} (\cite[Theorem 2.9]{2004-Chen}) Let $H \in \mathscr{A}$ and $\{g_n\}_{n=1}^{\infty} \subset \mathscr{H}$. Define $\tilde{f}^{(1)} = g_1$ and $\tilde{f}^{(n+1)} = H\tilde{f}^{(n)} + g_{n+1}$ for $n \ge 1$. If $g_n \uparrow g$ (respectively, $\sum_{n=1}^{\infty} g_n = g$), then $\tilde{f}^{(n)} \uparrow m_H g$ (respectively, $m_H g = \sum_{n=1}^{\infty} \tilde{f}^{(n)}$).
\end{mylem}

\begin{mylem}\label{Second-ite-geo} (\cite[Theorem 2.9]{2022-Wei}) Let $f^*$ be the minimal solution to \eqref{mini-solution}, and let $\hat{f}$ be a non-negative function satisfying 
\begin{equation*}
\hat{f}(x) \le (H\hat{f})(x) + g(x), \quad x \in \mathscr{X},
\end{equation*}
If $\hat{f} \le p f^*$ for some non-negative number $p$, then $\hat{f} \le f^*$.
\end{mylem}

\subsection{Lower bounds for moments of hitting times and sufficiency}
In this section, we prove the sufficiency of Theorem \ref{Non-1-ergo}, Theorem \ref{Non-stro}, and Proposition \ref{Dykin} using the minimal solution theory and Dynkin's formula.

For any set $A \in \mathscr{B}(\mathscr{X})$, define the first return time $\tau_A^+$ and first hitting time $\tau_A$ on $A$ as follows:
\[
\tau_A^+ := \inf\{n \ge 1 : X_n \in A\}, \qquad \tau_A := \inf\{n \ge 0 : X_n \in A\}.
\]
Denote
\[
L(x,A) = \sum_{n=1}^{\infty} \mathbb{P}_x(\tau_A^+ = n) = \mathbb{P}_x(\tau_A^+ < \infty).
\]
Before proceeding further, let us briefly outline the main ideas of our proofs.

Taking non-strong ergodicity as an example, it has an alternative characterization (\cite[Theorem 16.0.2]{2009-M-T}): for an aperiodic chain $\{X_n\}$, it is strongly ergodic if and only if there exists a petite set $C \in \mathscr{B}(\mathscr{X})$ such that $\sup_{x \in \mathscr{X}} \mathbb{E}_x \tau_C^+ < \infty$, and for every set $A \in \mathscr{B}^+(\mathscr{X})$, we have $\sup_{x \in \mathscr{X}} \mathbb{E}_x \tau_A^+ < \infty$. To verify that a chain is non-strongly ergodic, it suffices to show that $\sup_{x \in \mathscr{X}} \mathbb{E}_x \tau_A^+ = \infty$ for some set $A \in \mathscr{B}^+(\mathscr{X})$. We first derive a lower bound for the expectation of the return time using the minimal solution theory (see Proposition \ref{thLB}); an increasing sequence of such lower bounds then implies the desired result. On the other hand, a finite approximation method (see Proposition \ref{approx} and Lemmas \ref{Lem-nece-2}, \ref{lem-nece}) guarantees the existence of an increasing sequence of lower bounds, thus establishing the necessity of our conditions.

Additionally, by incorporating appropriate Lyapunov-like conditions, one can derive lower bounds for the hitting times of the Markov chain using Dynkin's formula. This allows us to obtain Lyapunov-like criteria for non-strong ergodicity.

  The following assertion, which relaxes the nonnegativity condition, is a modified version of \cite[Theorem 8.0.2]{2009-M-T}.

\begin{myprop}\label{transient-1} A $\psi$-irreducible Markov chain $\{X_n\}$ is transient if and only if there exist a set $C \in \mathscr{B}^+(\mathscr{X})$ and a function $V(x): \mathscr{X} \to \mathbb{R}$ with $\inf_{x \in \mathscr{X}} V(x) > -\infty$ such that
\[
PV(x) \le V(x) \quad \text{for all } x \in C^c,
\]
and
\begin{equation}\label{Trs1}
\left\{ x : V(x) < \inf_{y \in C} V(y) \right\} \in \mathscr{B}^+(\mathscr{X}).
\end{equation}
\end{myprop}

\textbf{Proof} The necessity follows directly from \cite[Theorem 8.0.2]{2009-M-T}. We now prove the sufficiency. By \cite[Theorem 8.3.6(ii)]{2009-M-T}, it suffices to show that there exist sets $C$ and $D \in \mathscr{B}^+(\mathscr{X})$ such that
\[
\mathbb{P}_x(\tau_C^+ < \infty) < 1 \quad \text{for all } x \in D.
\]
Without loss of generality, assume $V(x)$ is nonnegative (otherwise, replace $V(x)$ by $V(x) - \inf_{x \in \mathscr{X}} V(x)$). For the set $C \in \mathscr{B}^+(\mathscr{X})$, let $a = \inf_{x \in C} V(x)$ and define $D = \{ x : V(x) < a \}$. By \eqref{Trs1}, we have $a > 0$, $D \in \mathscr{B}^+(\mathscr{X})$, and $D \subseteq C^c$. For $x \in D$,
\[
\begin{aligned}
a > V(x) &\ge PV(x) = \int_{C} V(y) P(x, \mathrm{d}y) + \int_{C^c} V(y) P(x, \mathrm{d}y) \\
&\ge a P(x, C) + \int_{C^c} PV(y) P(x, \mathrm{d}y) \\
&\ge a P(x, C) + a \int_{C^c} P(y, C) P(x, \mathrm{d}y) + \int_{C^c} \int_{C^c} V(z_1) P(y, \mathrm{d}z_1) P(x, \mathrm{d}y) \\
&= a \mathbb{P}_x(\tau_C^+ = 1) + a \mathbb{P}_x(\tau_C^+ = 2) + \int_{C^c} \int_{C^c} V(z_1) P(y, \mathrm{d}z_1) P(x, \mathrm{d}y).
\end{aligned}
\]
By induction, for all $x \in D$ and $n \ge 1$, we obtain
\[
\begin{aligned}
a > V(x) > a \sum_{k=1}^n \mathbb{P}_x(\tau_C^+ = k) + \int_{C^c} \cdots \int_{C^c} V(z_{n-1}) P(z_{n-2}, \mathrm{d}z_{n-1}) \cdots P(x, \mathrm{d}y).
\end{aligned}
\]
Letting $n \to \infty$ yields
\[
\mathbb{P}_x(\tau_C^+ < \infty) < 1 \quad \text{for all } x \in D.
\]
$\Box$

 The following characterization of the moment of the return time via the minimal solution is essential for our application of the minimal solution theory.

\begin{myprop}\label{algbra-Moments} For each $A \in \mathscr{B}(\mathscr{X})$, $\mathbb{E}_x \tau_A^+ \mathbf{1}_{\{\tau_A^+ < \infty\}}$ is the minimal solution to the following equation:
\begin{equation}\label{am-3}
V(x) = \int_{A^c} V(y) P(x, \mathrm{d}y) + L(x, A), \qquad x \in \mathscr{X}.
\end{equation}
\end{myprop}

\textbf{Proof} We prove this using the second successive approximation scheme for minimal solutions. Set $V^{(1)}(x) = \mathbb{P}_x(\tau_A^+ = 1)$ and consider the inductive equation:
\[
V^{(n+1)}(x) = \int_{A^c} V^{(n)}(y) P(x, \mathrm{d}y) + \mathbb{P}_x(\tau_A^+ = n+1), \quad x \in \mathscr{X}.
\]
\iffalse\footnote{
Suppose that $V^{(k)}(x) = k\mathbb{P}_x(\tau_A^+ = k)$ holds for all $k \le n$. Then for $k = n+1$:
\[
\begin{aligned}
V^{(n+1)}(x) &= \int_{A^c} V^{(n)}(y) P(x, \mathrm{d}y) + \mathbb{P}_x(\tau_A^+ = n+1) \\
&= \int_{A^c} n\mathbb{P}_y(\tau_A^+ = n) P(x, \mathrm{d}y) + \mathbb{P}_x(\tau_A^+ = n+1) \\
&= (n+1)\mathbb{P}_x(\tau_A^+ = n+1).
\end{aligned}
\]} \fi
This equation has a solution $V^{(n)}(x) = n\mathbb{P}_x(\tau_A^+ = n)$ for all $n \ge 1$. Note that $L(x, A) = \sum_{n=1}^{\infty} \mathbb{P}_x(\tau_A^+ = n)$. By Lemma \ref{Second-ite},
\[
V^*(x) := \sum_{n=1}^{\infty} V^{(n)}(x) = \mathbb{E}_x \tau_A^+ \mathbf{1}_{\{\tau_A^+ < \infty\}}, \quad x \in \mathscr{X},
\]
is the minimal nonnegative solution to \eqref{am-3}. $\Box$

Recall that $\{X_n\}$ is recurrent if and only if $L(x, A) = 1$ for all $x \in \mathscr{X}$ and $A \in \mathscr{B}^+(\mathscr{X})$ (\cite[Theorem 8.3.6, Propositions 9.1.1, 9.1.4, 9.1.7]{2009-M-T}). The following assertion provides a "lower bound" for the moment of the return time.

\begin{myprop}\label{thLB} Assume that the Markov chain $\{X_n\}$ with one-step transition probability kernel $\{P(x, \cdot) : x \in \mathscr{X}\}$ is $\psi$-irreducible and recurrent. Let $A \in \mathscr{B}^+(\mathscr{X})$ and let $V(x): \mathscr{X} \to \mathbb{R}$ be a function with $\sup_{x \in A^c} V(x) < \infty$ that satisfies
\begin{equation}\label{thLB-1-0}
V(x) \le \int_{A^c} V(y) P(x, \mathrm{d}y) + 1.
\end{equation}
Then the following hold:
\begin{itemize}
\item[(1)] If \eqref{thLB-1-0} holds for all $x \in A^c$, then $V(x) \le \mathbb{E}_x \tau_A^+ + \sup_{x \in A} V(x)$ for all $x \in \mathscr{X}$. Moreover, if $V(x) = 0$ for all $x \in A$, then
\[
V(x) \le \mathbb{E}_x \tau_A^+ \quad \text{for all } x \in \mathscr{X}.
\]
\item[(2)] If \eqref{thLB-1-0} holds for all $x \in A^c$, then
\[
\psi\left( \left\{ x \in A^c : V(x) > \mathbb{E}_x \tau_A^+ \right\} \right) = 0.
\]
\item[(3)] If \eqref{thLB-1-0} holds for all $x \in \mathscr{X}$, then
\[
\psi\left( \left\{ y \in A : V(y) > \mathbb{E}_y \tau_A^+ \right\} \right) = 0.
\]
Moreover,
\[
\psi\left( \left\{ y \in \mathscr{X} : V(y) > \mathbb{E}_y \tau_A^+ \right\} \right) = 0.
\]
\end{itemize}
\end{myprop}\textbf{Proof} (1) Since
\[
PV(x) = P(V\mathbf{1}_A)(x) + PV\mathbf{1}_{A^c}(x) \ge P(V\mathbf{1}_A)(x),
\]
and \eqref{thLB-1-0} holds for all $x \in A^c$, by Dynkin's formula and the Markov property (\cite[Theorem 11.3.1]{2009-M-T}), for all $x \in A^c$ and $m \ge 1$,
\begin{equation}\label{d-k-0}
\begin{aligned}
\mathbb{E}_x V(X_{m \wedge \tau_A^+})
&= V(x) + \mathbb{E}_x \bigg[ \sum_{i=0}^{m \wedge \tau_A^+ - 1} \bigg( PV(X_i) - V(X_i) \bigg) \bigg] \\
&\ge V(x) - \mathbb{E}_x \big( m \wedge \tau_A^+ \big).
\end{aligned}
\end{equation}

Since $\{X_n\}$ is Harris recurrent, $\mathbb{P}_x(\tau_A < \infty) = 1$ for all $A \in \mathscr{B}^+(\mathscr{X})$ and $x \in \mathscr{X}$. Thus,
\[
\mathbb{P}_x(\tau_A > m) \to 0 \quad \text{as } m \to \infty \quad \text{for all } x \in \mathscr{X}.
\]

Note that $\mathbb{E}_x V(X_{m \wedge \tau_A^+}) = \mathbb{E}_x V(X_m) \mathbf{1}_{\{m < \tau_A^+\}} + \mathbb{E}_x V(X_{\tau_A^+}) \mathbf{1}_{\{m \ge \tau_A^+\}}$. Since $\sup_{x \in A^c} V(x) < \infty$, we have
\[
\varlimsup_{m \to \infty} \mathbb{E}_x V(X_{m \wedge \tau_A^+}) \le \mathbb{E}_x V(X_{\tau_A^+}) \le \sup_{x \in A} V(x) \quad \text{for all } x \in \mathscr{X}.
\]

Moreover, by \eqref{d-k-0}, $\mathbb{E}_x \tau_A^+ + \sup_{x \in A} V(x) \ge V(x)$ for all $x \in \mathscr{X}$, so the required assertions hold.

(2) Since $\{X_n\}$ is recurrent, $L(x, A) = 1$ for all $x \in \mathscr{X}$. By Proposition \ref{algbra-Moments},
\begin{equation}\label{thLB-2-0}
\mathbb{E}_x \tau_A^+ = \int_{A^c} \mathbb{E}_y \tau_A^+ \, P(x, \mathrm{d}y) + 1 \quad \text{for all } x \in \mathscr{X}.
\end{equation}

Let
\[
Z(x) = \big[ \mathbb{E}_x \tau_A^+ - V(x) \big] \mathbf{1}_{A^c}(x).
\]

Note that $Z(x)$ is well-defined since $\sup_{x \in A^c} V(x) < \infty$, which also implies $\inf_{x \in \mathscr{X}} Z(x) > -\infty$. For all $x \in A^c$, by \eqref{thLB-1-0} and \eqref{thLB-2-0},
\[
\begin{aligned}
PZ(x) &= \int_{A^c} Z(y) P(x, \mathrm{d}y) \\
&= \int_{A^c} \mathbb{E}_y \tau_A^+ \, P(x, \mathrm{d}y) - \int_{A^c} V(y) P(x, \mathrm{d}y) \\
&\le \mathbb{E}_x \tau_A^+ - 1 - \big[ V(x) - 1 \big],
\end{aligned}
\]
so we obtain $PZ(x) \le Z(x)$ for all $x \in A^c$. Since $\{X_n\}$ is recurrent, Proposition \ref{transient-1} implies
\[
\psi\left( \left\{ x : Z(x) < \inf_{y \in A} Z(y) \right\} \right) = 0.
\]

That is, $\psi(D) = 0$ where
\begin{equation}\label{set-D}
D = \left\{ y \in A^c : V(y) > \mathbb{E}_y \tau_A^+ \right\}.
\end{equation}

(3) Let $\pi$ denote the invariant measure of $\{X_n\}$ satisfying $\pi P = \pi$. By \cite[Theorem 10.4.9]{2009-M-T}, $\pi \ll \psi$ and $\psi \ll \pi$, so $P(x, D) = 0$ for $\pi$-almost all $x$, where $D$ is defined in \eqref{set-D}. There exists a measurable set $B$ with $\pi(B) = 0$ (equivalently, $\psi(B) = 0$) such that for all $x \in A \setminus B$, $P(x, D) = 0$ and
\[
\begin{aligned}
V(x) &\le \int_{A^c} V(y) P(x, \mathrm{d}y) + 1 \quad (\text{since \eqref{thLB-1-0} holds for } x \in A) \\
&\le \int_{A^c \setminus D} V(y) P(x, \mathrm{d}y) + 1 + \int_{D} V(y) P(x, \mathrm{d}y) \\
&\le \int_{A^c} \mathbb{E}_y \tau_A^+ P(x, \mathrm{d}y) + 1 + \int_{D} \bigg( V(y) - \mathbb{E}_y \tau_A^+ \bigg) P(x, \mathrm{d}y) \quad (\text{by definition of } D) \\
&= \mathbb{E}_x \tau_A^+ \quad \text{for all } x \in A \setminus B.
\end{aligned}
\]
The last equality follows from Proposition \ref{algbra-Moments} and the fact that $\sup_{y \in D} V(y) < +\infty$. Hence,
\[
\psi\left( \left\{ x \in A : V(x) > \mathbb{E}_x \tau_A^+ \right\} \right) = 0.
\]

\textbf{Proof of sufficiency of Theorem \ref{Non-1-ergo}} Suppose $\{X_n\}$ is ergodic. Then $\{X_n\}$ is recurrent. By parts (1), (2) of Proposition \ref{thLB}, we have
\[
\psi(N_n) = 0 \quad \text{for all } n \ge 1,
\]
where $N_n := \left\{ y \in \mathscr{X} : V^{(n)}(y) > \mathbb{E}_y \tau_A^+ \right\}$. Let $N = \cup_{n=1}^{\infty} N_n$. Then $\psi(N) = 0$ (equivalently, $\mu(N) = 0$ by \cite[Theorem 10.4.9]{2009-M-T}, where $\mu$ is the invariant measure) and
\[
\mathbb{E}_y \tau_A^+ \ge V^{(n)}(y) \quad \text{for all } y \in \mathscr{X} \setminus N \text{ and } n \ge 1.
\]
Thus,
\[
\int_{A} \mathbb{E}_y \tau_A^+ \mu(\mathrm{d}y) \ge \int_{A} \sup_{n \ge 1} V^{(n)}(y) \mu(\mathrm{d}y) \ge \sup_{n \ge 1} \int_A V^{(n)}(y) \mu(\mathrm{d}y).
\]
Since $\pi$ and $\mu$ are equivalent, condition (3) of Theorem \ref{Non-1-ergo} implies
\[
\int_{A} \mathbb{E}_y \tau_A^+ \mu(\mathrm{d}y) = \infty.
\]
This contradicts \cite[Theorem 10.4.10 (i)]{2009-M-T}, so the required assertion holds. $\Box$

\begin{myprop}\label{complement} Let $\{X_n\}$ be a $\psi$-irreducible Markov chain. If there exist sets $A, B \in \mathscr{B}^+(\mathscr{X})$ such that
\begin{equation}\label{comp}
\left\{ x \in A : \mathbb{E}_x \tau_B^+ = \infty \right\} \in \mathscr{B}^+(\mathscr{X}),
\end{equation}
then $\{X_n\}$ is non-ergodic.
\end{myprop}

\textbf{Proof} Suppose $\{X_n\}$ is ergodic. By \cite[Theorem 11.1.4]{2009-M-T}, there exists a full set $S$ (i.e., $\psi(S^c) = 0$) and an absorbing set (i.e., $P(x, S) = 1$ for all $x \in S$) that admits a countable cover $\{S_n\}$ such that
\[
\sup_{x \in S_n} \mathbb{E}_x \tau_B^+ < \infty \quad \text{for all } B \in \mathscr{B}^+(\mathscr{X}) \text{ and } n \ge 1.
\]
Therefore,
\[
\psi\left( \left\{ x \in \mathscr{X} : \mathbb{E}_x \tau_B^+ = \infty \right\} \right) = \psi\left( \left\{ x \in S : \mathbb{E}_x \tau_B^+ = \infty \right\} \cup S^c \right) = 0 \quad \text{for all } B \in \mathscr{B}^+(\mathscr{X}),
\]
which contradicts \eqref{comp}. $\Box$

\textbf{Proof of Corollary \ref{Non-1-ergo-cor}} By Proposition \ref{thLB}, if \eqref{thLB-1-0} holds for all $x \in \mathscr{X}$, there exists a set $N$ with $\psi(N) = 0$ such that
\[
V^{(n)}(x) \le \mathbb{E}_x \tau_A^+ \quad \text{for all } x \in \mathscr{X} \setminus N \text{ and } n \ge 1.
\]
Hence, for any $B \in \mathscr{B}(\mathscr{X})$,
\[
\left\{ x \in B : \sup_{n \ge 1} V^{(n)}(x) = \infty \right\} \subseteq \left\{ x \in B : \mathbb{E}_x \tau_A^+ = \infty \right\}.
\]
If \eqref{thLB-1-0} holds only for $x \in A^c$, we have $N \subseteq A^c$. The first required assertion follows immediately from condition (3) of the corollary and Proposition \ref{complement}. $\Box$

\textbf{Proof of sufficiency of Theorem \ref{Non-stro}} Suppose $\{X_n\}$ is strongly ergodic. Then $\mathbb{P}_x(\tau_A^+ < \infty) = 1$ for all $A \in \mathscr{B}^+(\mathscr{X})$ and $x \in \mathscr{X}$. By Proposition \ref{thLB} and conditions (1), (2) of Theorem \ref{Non-stro},
\[
\mathbb{E}_x \tau_A^+ \ge V^{(n)}(x) \quad \text{for all } x \in \mathscr{X} \text{ and } n \ge 1.
\]
Hence,
\[
\sup_{x \in A^c} \mathbb{E}_x \tau_A^+ \ge \sup_{\substack{x \in A^c \\ n \ge 1}} V^{(n)}(x) = \infty
\]
by condition (3) of Theorem \ref{Non-stro}. This contradicts the fact that $\sup_{x \in \mathscr{X}} \mathbb{E}_x \tau_A^+ < \infty$ for strongly ergodic chains. Therefore, $\{X_n\}$ is non-strongly ergodic. $\Box$

\textbf{Proof of Proposition \ref{Dykin}} Using an argument similar to that in the proof of Theorem \ref{Non-stro}, for all $x \in E_n \setminus A$ and $m \ge 1$, we have
\[
\mathbb{E}_x V(X_{m \wedge \tau_A \wedge \tau_{E_n^c}}) \ge V(x) - \mathbb{E}_x \big( m \wedge \tau_A \wedge \tau_{E_n^c} \big),
\]
and
\begin{equation}\label{fz1}
\mathbb{E}_x W(X_{m \wedge \tau_A \wedge \tau_{E_n^c}}) \le W(x).
\end{equation}

Suppose the chain is strongly ergodic. Then $\mathbb{P}_x(\tau_A < \infty) = 1$ for all $A \in \mathscr{B}^+(\mathscr{X})$, and $\lim_{m \to \infty} \mathbb{P}_x(\tau_A > m) = 0$ for all $x \in \mathscr{X}$. By \eqref{fz1} and Fatou's lemma,
\begin{equation}\label{fz2}
\mathbb{E}_x W(X_{\tau_A \wedge \tau_{E_n^c}}) \le \varliminf_{m \to \infty} \mathbb{E}_x W(X_{m \wedge \tau_A \wedge \tau_{E_n^c}}) \le W(x) \quad \text{for all } x \in E_n \setminus A.
\end{equation}

By the local boundedness of $V$, we have
\[
|V(X_m)| \mathbf{1}_{\{m < \tau_A \wedge \tau_{E_n^c}\}} \le \sup_{x \in E_n \setminus A} |V(x)| < \infty.
\]
Thus, $\varlimsup_{m \to \infty} \mathbb{E}_x V(X_m) \mathbf{1}_{\{m < \tau_A \wedge \tau_{E_n^c}\}} = 0$ for all $x \in E_n \setminus A$, and
\[
\varlimsup_{m \to \infty} \mathbb{E}_x V(X_{m \wedge \tau_A \wedge \tau_{E_n^c}}) = \varlimsup_{m \to \infty} \mathbb{E}_x V(X_{\tau_A \wedge \tau_{E_n^c}}) \mathbf{1}_{\{m \ge \tau_A \wedge \tau_{E_n^c}\}} \le \mathbb{E}_x V(X_{\tau_A \wedge \tau_{E_n^c}}).
\]

Therefore,
\begin{equation}\label{js}
\mathbb{E}_x (\tau_A \wedge \tau_{E_n^c}) \ge V(x) - \mathbb{E}_x V(X_{\tau_A \wedge \tau_{E_n^c}}).
\end{equation}

Since $W$ is nonnegative, \eqref{fz2} implies
\[
\mathbb{E}_x \bigg[ W(X_{\tau_{E_n^c}}) \mathbf{1}_{\{\tau_A \ge \tau_{E_n^c}\}} \bigg] \le W(x) \quad \text{for all } x \in E_n \setminus A.
\]

Furthermore,
\[
\begin{aligned}
\mathbb{E}_x (\tau_A \wedge \tau_{E_n^c}) &\ge V(x) - \mathbb{E}_x V(X_{\tau_A \wedge \tau_{E_n^c}}) \\
&\ge V(x) - \mathbb{E}_x V(X_{\tau_A}) \mathbf{1}_{\{\tau_A < \tau_{E_n^c}\}} - \bigg( \sup_{y \in E_n^c} \frac{V(y)}{W(y)} \bigg) \mathbb{E}_x \big[ W(X_{\tau_{E_n^c}}) \mathbf{1}_{\{\tau_A \ge \tau_{E_n^c}\}} \big] \\
&\ge V(x) - \mathbb{E}_x V(X_{\tau_A}) \mathbf{1}_{\{\tau_A < \tau_{E_n^c}\}} - \bigg( \sup_{y \in E_n^c} \frac{V(y)}{W(y)} \bigg) W(x) \quad \text{for all } x \in E_n \setminus A.
\end{aligned}
\]

As $n \to \infty$, $\tau_{E_n^c} \uparrow \infty$. By condition (1) of the proposition and the fact that $\sup_{x \in A} V(x) < \infty$, we obtain
\[
\sup_{x \in A^c} \mathbb{E}_x \tau_A \ge \sup_{x \in A^c} V(x) - \sup_{x \in A} V(x) = \infty.
\]

This contradicts the property of strongly ergodic chains that $\sup_{x \in \mathscr{X}} \mathbb{E}_x \tau_A^+ < \infty$. Therefore, $\{X_n\}$ is non-strongly ergodic. $\Box$

\iffalse
\begin{myrem}
\begin{itemize}
\item[1)] Since $\sup_{y \in E_n \setminus A} V(y) \to \infty$ as $n \to \infty$, we analyze $\mathbb{E}_x V(X_{\tau_A \wedge \tau_{E_n^c}})$ in \eqref{js} by introducing an additional test function $W$.
\item[2)] The boundedness of $V(X_m) \mathbf{1}_{\{m < \tau_A\}}$ (i.e., $V(X_m) \mathbf{1}_{\{m < \tau_A\}} \le \sup_{x \in A^c} V(x)$) cannot be guaranteed. Thus, the dominated convergence theorem is inapplicable if only $V$ is considered.
\end{itemize}
\end{myrem}
\fi

\subsection{Approximation for the moments of hitting times; necessity of Theorems \ref{Non-1-ergo} \ref{Non-stro} and proof of Theorem \ref{Non-geo}}

\subsubsection{Approximation of Markov chain}
\begin{myprop}\label{approx} Assume that the Markov chain $\{X_n\}$ with Markov semigroup $\{P(x, \cdot)\}$ on a locally compact separable metric space $(\mathscr{X}, \mathscr{B}(\mathscr{X}))$ is $\psi$-irreducible, aperiodic, Feller, and $\mathrm{supp}\,\psi$ has nonempty interior. Let $E \subseteq \mathscr{X}$ be a compact set with $0 < \psi(E) < \infty$. For any $A \in \mathscr{B}^+(\mathscr{X})$ with $A \subseteq E$ being a compact set, define $\{\hat{P}(x, \cdot)\}$ on $(E, \mathscr{B}(E))$ as follows:
\begin{equation}\label{constr}
\hat{P}(x, B) = \bigg[ P(x, B) + P(x, E^c) \frac{\psi(B \cap A)}{\psi(A)} \bigg] \mathbf{1}_{E \setminus A}(x) + \frac{\psi(B)}{\psi(E)} \mathbf{1}_A(x)
\end{equation}
for all $B \in \mathscr{B}(E)$. Then $\{\hat{P}(x, \cdot) : x \in E\}$ is a transition probability, and the Markov chain generated by $\{\hat{P}(x, \cdot) : x \in E\}$ is $\psi$-irreducible, aperiodic, and strongly ergodic.
\end{myprop}

\textbf{Proof} It is easy to verify that $\hat{P}(x, E) = 1$ and $\{\hat{P}(x, B) : x \in E, B \in \mathscr{B}(E)\}$ is a transition probability on $(E, \mathscr{B}(E))$. Let $\{\hat{X}_n\}$ denote the Markov chain on the state space $(E, \mathscr{B}(E))$ with one-step transition probability $\{\hat{P}(x, B) : x \in E, B \in \mathscr{B}(E)\}$. For $n \ge 2$ and $B \in \mathscr{B}(E)$, let $\hat{P}_n(x, B) = \mathbb{P}_x(\hat{X}_n \in B)$, with $\hat{P}_1 = \hat{P}$.

We first prove that $\{\hat{X}_n\}$ is $\psi$-irreducible, i.e., for all $x \in E$ and $B \in \mathscr{B}(E)$ with $\psi(B) > 0$, there exists $n \ge 1$ such that $\hat{P}_n(x, B) > 0$ (by \cite[Proposition 4.2.1]{2009-M-T}).

1) For all $x \in A$ and $B \in \mathscr{B}(E)$ with $\psi(B) > 0$, it is obvious that $\hat{P}(x, B) > 0$.

2) We prove that $\sum_{n=1}^{\infty} \hat{P}_n(x, A) > 0$ for all $x \in E \cap A^c$ (equivalently, for all $x \in E \cap A^c$, there exists $n \ge 1$ such that $\hat{P}_n(x, A) > 0$). Suppose this is false, so there exists $x_1 \in E \cap A^c$ such that $\hat{P}_n(x_1, A) = 0$ for all $n \ge 1$. We derive a contradiction as follows:

By \eqref{constr}, $\hat{P}_1(x_1, A) = \hat{P}(x_1, A) = P(x_1, A) + P(x_1, E^c) = 0$, which implies $P(x_1, A) = 0$ and $P(x_1, E^c) = 0$. Moreover, $P(x_1, A^c \cap E) = 1$ and $\hat{P}(x_1, A^c \cap E) = 1$, with
\[
\hat{P}(x_1, A^c \cap E \cap C) = P(x_1, A^c \cap E \cap C) \quad \text{for all } C \in \mathscr{B}(\mathscr{X}).
\]

By the Chapman-Kolmogorov equation, for all $C \in \mathscr{B}(\mathscr{X})$,
\[
\begin{aligned}
\hat{P}_2(x_1, A^c \cap E \cap C) &= \int_{E \cap A^c} \hat{P}(x_1, \mathrm{d}y) \hat{P}(y, A^c \cap E \cap C) \quad (\text{since } \hat{P}(x_1, A^c \cap E) = 1) \\
&= \int_{E \cap A^c} P(x_1, \mathrm{d}y) P(y, A^c \cap E \cap C) \quad (\text{by } \eqref{constr}) \\
&= P_2(x_1, E \cap A^c \cap C) \quad (\text{since } P(x_1, A^c \cap E) = 1).
\end{aligned}
\]

Similarly, $\hat{P}_2(x_1, E^c) = 0$ by $P(x_1, E^c) = 0$ and the Chapman-Kolmogorov equation. Since $\hat{P}_2(x_1, A) = 0$, we have $\hat{P}_2(x_1, A^c \cap E) = P_2(x_1, E \cap A^c) = 1$.

Inductively, assume that for all $n \le k$,
\begin{equation}\label{induc1}
\hat{P}_n(x_1, A^c \cap E \cap C) = P_n(x_1, A^c \cap E \cap C) \quad \text{for all } C \in \mathscr{B}(\mathscr{X})
\end{equation}
and $\hat{P}_n(x_1, A^c \cap E) = P_n(x_1, E \cap A^c) = 1$. Then
\[
\begin{aligned}
\hat{P}_{k+1}(x_1, A^c \cap E \cap C) &= \int_E \hat{P}_k(x_1, \mathrm{d}y) \hat{P}(y, A^c \cap E \cap C) \\
&= \int_{E \cap A^c} \hat{P}_k(x_1, \mathrm{d}y) \hat{P}(y, A^c \cap E \cap C) \quad (\text{since } \hat{P}_k(x_1, A) = 0) \\
&= \int_{E \cap A^c} P_k(x_1, \mathrm{d}y) P(y, A^c \cap E \cap C) \quad (\text{by } \eqref{constr} \text{ and } \eqref{induc1}) \\
&= P_{k+1}(x_1, E \cap A^c \cap C) \quad (\text{since } P_k(x_1, A^c \cap E) = 1).
\end{aligned}
\]

Thus, if there exists $x_1 \in E \cap A^c$ with $\hat{P}_n(x_1, A) = 0$ for all $n \ge 1$, then $\hat{P}_n(x_1, A^c \cap E) = P_n(x_1, E \cap A^c) = 1$ for all $n \ge 1$. This implies $\sum_{n=1}^{\infty} P_n(x_1, E^c \cup A) = 0$, contradicting the irreducibility of $\{X_n\}$.

3) We prove that $\sum_{n=1}^{\infty} \hat{P}_n(x, B) > 0$ for all $x \in E \cap A^c$ and $B \in \mathscr{B}(E)$ with $\psi(B) > 0$. By part 2), there exists $n \ge 1$ such that $\hat{P}_n(x, A) > 0$. Moreover,
\[
\hat{P}_{n+1}(x, B) \ge \int_A \hat{P}_n(x, \mathrm{d}y) \hat{P}(y, B) = \frac{\psi(B)}{\psi(E)} \hat{P}_n(x, A) > 0.
\]
Thus, $\{\hat{X}_n\}$ is $\psi$-irreducible.

Since $A \subseteq E$ is compact, $\{\hat{X}_n\}$ is weakly Feller by definition. By \cite[Theorems 6.2.9, 16.2.5]{2009-M-T}, it suffices to prove aperiodicity. Let $\hat{\nu}_1(\cdot) = \psi(\cdot)/\psi(E)$ be a measure. For all $x \in A$ and $B \in \mathscr{B}(E)$, $\hat{P}(x, B) = \hat{\nu}_1(B)$, so $A$ is $\hat{\nu}_1$-small by definition (see \cite[Section 5.2]{2009-M-T}). Since $\hat{\nu}_1(A) > 0$, $\{\hat{X}_n\}$ is strongly aperiodic by definition (see \cite[Section 5.4.3]{2009-M-T}). Hence, the required assertion holds. $\Box$

\begin{myexam}
Consider the Markov chain $X_{n+1} = aX_n + W_{n+1}$ on $\mathbb{R}$ with $X_0 = 0$, $|a| < 1$, and $a \neq 0$. For any $m \ge 1$, let $E_m = \left[ -m + \frac{1}{m}, m - \frac{1}{m} \right]$, so $E_m \uparrow \mathbb{R}$ as $m \to \infty$. Its one-step transition kernel $\{P(x, \cdot) : x \in \mathbb{R}\}$ has the intensity function
\[
p(x, y) = \frac{1}{\sqrt{2\pi}} \exp\left\{ -\frac{(y - ax)^2}{2} \right\}.
\]
Let $A = E_1$, which satisfies $\mu_{\text{Leb}}(A) > 0$ (where $\mu_{\text{Leb}}$ denotes the Lebesgue measure). Define
\begin{equation}\label{LZZ1}
p^{(m)}(x, y) = 
\begin{cases}
p(x, y) + P(x, E_m^c) \mathbf{1}_A(y) \frac{y}{\mu_{\text{Leb}}(A)}, & x \in E_m \setminus A; \\
\frac{1}{\mu_{\text{Leb}}(E_m)}, & x \in A;
\end{cases}
\end{equation}
and let $P^{(m)}(x, \mathrm{d}y) = p^{(m)}(x, y) \mathrm{d}y$. Then $\{P^{(m)}(x, \mathrm{d}y)\}$ is a transition kernel on $E_m$. Let $\{X_n^{(m)}\}$ be the Markov chain generated by $\{P^{(m)}(x, \cdot)\}$. It is easy to verify that $\{P^{(m)}(x, \mathrm{d}y)\}$ is irreducible with respect to the Lebesgue measure on $E_m$. Furthermore, $\{X_n^{(m)}\}$ is strongly ergodic by Proposition \ref{approx}.
\end{myexam}

\begin{myrem}
It is easy to see that $\{p^{(m)}(x, y)\}$ in \eqref{LZZ1} is increasing in $m$ and converges to $p(x, y)$ as $m \to \infty$.
\end{myrem}

\subsubsection{Necessity of Theorems \ref{Non-1-ergo} and \ref{Non-stro}}
In this subsection, we assume the locally compact separable metric space $\mathscr{X} = \cup_{m=1}^{\infty} E_m$, where $\{E_m\}_{m \ge 1}$ is an increasing sequence of compact sets. The Markov chain $\{X_n\}$ with Markov semigroup $\{P(x, \cdot)\}$ on $(\mathscr{X}, \mathscr{B}(\mathscr{X}))$ always satisfies Assumption 1 (i.e., $\psi$-irreducible, aperiodic, Feller, and $\mathrm{supp}\,\psi$ has nonempty interior).

The following lemma is an immediate consequence of Proposition \ref{approx}.
\begin{mylem}\label{Lem-nece-2}
Let $A \in \mathscr{B}^+(E_1)$ be a compact set. For each $m \ge 1$, let $\{X_n^{(m)} : n \ge 0\}$ be a Markov chain on the state space $E_m$ with one-step transition semigroup $P^{(m)}(x, \cdot)$, defined as $\hat{P}$ in \eqref{constr} by replacing $E$ with $E_m$. Then the Markov chain $\{X_n^{(m)}\}$ on $(E_m, \mathscr{B}(E_m))$ is strongly ergodic.
\end{mylem}

For a set $B \in \mathscr{B}(E_m)$, let $\tau^{(m)+}_B = \inf\{n \ge 1 : X_n^{(m)} \in B\}$ for the Markov chain $\{X_n^{(m)}\}$ in Lemma \ref{Lem-nece-2}. Note that $\mathbb{P}_x(\tau^{(m)+}_B < \infty) = 1$ for all $x \in E_m$. Consider the equation
\begin{equation}\label{truc}
V(x) = \int_{B^c \cap E_m} V(y) P^{(m)}(x, \mathrm{d}y) + 1, \quad x \in E_m.
\end{equation}
By Proposition \ref{algbra-Moments}, the minimal solution to \eqref{truc} is $\mathbb{E}_x^{(m)} \tau^{(m)+}_B$ for all $x \in E_m$, where $\mathbb{E}_x^{(m)}$ denotes the expectation with respect to the law of $\{X_n^{(m)}\}$ starting at $x$.

By Proposition \ref{approx}, we immediately obtain the following lemma:
\begin{mylem}\label{lem-nece}
Let $A \in \mathscr{B}^+(E_1)$ be a compact set and assume $\{X_n\}$ is recurrent. Then the following assertions hold:
\begin{itemize}
\item[(1)] $\sup_{x \in A^c \cap E_m} \mathbb{E}_x^{(m)} \tau^{(m)+}_A < \infty$ for all $m \ge 1$;
\item[(2)] $\mathbb{E}_x \tau_A^+ = \lim_{m \to \infty} \uparrow \mathbb{E}_x^{(m)} \tau^{(m)+}_A \mathbf{1}_{E_m}(x)$ for all $x \in \mathscr{X}$, and the sequence $\left\{ \sup_{x \in A^c \cap E_m} \mathbb{E}_x^{(m)} \tau^{(m)+}_A \right\}_{m \ge 1}$ is increasing;
\item[(3)] $\left\{ \sup_{x \in A^c \cap E_m} \mathbb{E}_x^{(m)} \tau^{(m)+}_A \right\}_{m \ge 1}$ is bounded if and only if $\{\mathbb{E}_x \tau_A^+ : x \in A^c\}$ is bounded. That is, $\sup_{x \in A^c} \mathbb{E}_x \tau_A^+ < \infty$ if and only if $\sup_{m \ge 1} \sup_{x \in A^c \cap E_m} \mathbb{E}_x^{(m)} \tau^{(m)+}_A < \infty$.
\end{itemize}
\end{mylem}

\textbf{Proof} (1) By Lemma \ref{Lem-nece-2}, $\{X_n^{(m)}\}$ is strongly ergodic. Since $A \in \mathscr{B}^+(E_m)$, we have $L^{(m)}(x, A) = 1$ (where $L^{(m)}(x, A) = \mathbb{P}_x(\tau^{(m)+}_A < \infty)$) and $\sup_{x \in A^c \cap E_m} \mathbb{E}_x^{(m)} \tau^{(m)+}_A < \infty$.
(2) Since $A \subset E_1$ and $\{X_n^{(m)}\}$ is Harris recurrent, by Proposition \ref{algbra-Moments}, $\mathbb{E}_x^{(m)} \tau^{(m)+}_A$ for all $x \in E_m$ is the minimal solution to
\[
V(x) = \int_{E_m \cap A^c} V(y) P^{(m)}(x, \mathrm{d}y) + 1, \qquad x \in E_m.
\]
Equivalently, $\widetilde{V}_A^{*(m)}(x) := \mathbb{E}_x^{(m)} \tau^{(m)+}_A \mathbf{1}_{E_m}(x)$ is the minimal solution to
\[
V(x) = \int_{E_m \cap A^c} V(y) P^{(m)}(x, \mathrm{d}y) \mathbf{1}_{E_m}(x) + \mathbf{1}_{E_m}(x), \qquad x \in \mathscr{X}.
\]
Since $E_m \uparrow \mathscr{X}$ as $m \to \infty$, we have
\[
\int_{A^c} V(y) \mathbf{1}_{E_m}(y) P^{(m)}(x, \mathrm{d}y) \mathbf{1}_{E_m}(x) \uparrow \int_{A^c} V(y) P(x, \mathrm{d}y) \quad \text{for all } x \in \mathscr{X}.
\]
Note that $\mathbb{E}_x \tau_A^+$ is the minimal solution to
\[
V(x) = \int_{A^c} V(y) P(x, \mathrm{d}y) + 1, \qquad x \in \mathscr{X}.
\]
By Lemma \ref{Monot}, $\lim_{m \to \infty} \uparrow \widetilde{V}_A^{*(m)}(x) = \mathbb{E}_x \tau_A^+$ for all $x \in \mathscr{X}$. Thus, $\mathbb{E}_x^{(m)} \tau^{(m)+}_A \mathbf{1}_{E_m}(x) \uparrow \mathbb{E}_x \tau_A^+$ for all $x \in \mathscr{X}$, and $\sup_{x \in A^c \cap E_m} \mathbb{E}_x^{(m)} \tau^{(m)+}_A$ is increasing in $m$.

(3) By part (2), the assertion holds trivially. $\Box$

\textbf{Proof of necessity of Theorem \ref{Non-1-ergo}} Suppose $\{X_n\}$ is non-ergodic but recurrent, with invariant measure $\mu$. By \cite[Theorem 10.4.10]{2009-M-T}, for every set $A$ with $\mu(A) > 0$,
\[
\int_A \mathbb{E}_y \tau_A^+ \mu(\mathrm{d}y) = \infty.
\]
Without loss of generality, choose a compact set $A \subseteq E_1$ with $\mu(A) > 0$ (this is feasible because $\mathscr{X}$ is a locally compact separable space). Define
\[
k_m(x) = \int_{E_m \cap A^c} \mathbb{E}_y^{(m)} \tau^{(m)+}_A P(x, \mathrm{d}y) + 1
\]
and
\[
W^{(m)}(x) = 
\begin{cases}
\mathbb{E}_x^{(m)} \tau^{(m)+}_A, & x \in E_m \cap A^c, \\
k_m(x), & x \notin E_m \cap A^c.
\end{cases}
\]

By Lemma \ref{lem-nece}, $\sup_{x \in A^c} k_m(x) < \infty$ for all $m \ge 1$, so condition (1) of Theorem \ref{Non-1-ergo} holds. Next, we compute
\[
\int_{A^c} W^{(m)}(y) P(x, \mathrm{d}y) = \int_{E_m \cap A^c} \mathbb{E}_y^{(m)} \tau^{(m)+}_A P(x, \mathrm{d}y) + \int_{E_m^c \cup A} k_m(y) P(x, \mathrm{d}y).
\]

For $x \in E_m \cap A^c$, $P^{(m)}(x, \mathrm{d}y) = P(x, \mathrm{d}y)$, and by Proposition \ref{algbra-Moments},
\[
\int_{E_m \cap A^c} \mathbb{E}_y^{(m)} \tau^{(m)+}_A P^{(m)}(x, \mathrm{d}y) = \mathbb{E}_x^{(m)} \tau^{(m)+}_A - 1.
\]
Thus,
\[
\int_{A^c} W^{(m)}(y) P(x, \mathrm{d}y) \ge W^{(m)}(x) - 1.
\]

For $x \notin E_m \cap A^c$,
\[
\int_{A^c} W^{(m)}(y) P(x, \mathrm{d}y) = k_m(x) - 1 + \int_{E_m^c \cup A} k_m(y) P(x, \mathrm{d}y) \ge W^{(m)}(x) - 1,
\]
so condition (2) of Theorem \ref{Non-1-ergo} holds.

It is straightforward to verify that
\[
\sup_{m \ge 1} W^{(m)}(x) = \mathbb{E}_x \tau_A^+ \quad \text{for all } x \in A^c.
\]
Moreover,
\[
\int_A \sup_{m \ge 1} W^{(m)}(y) \psi(\mathrm{d}y) = \int_A \mathbb{E}_y \tau_A^+ \psi(\mathrm{d}y) = \infty.
\]
The last equality follows from $\int_A \mathbb{E}_y \tau_A^+ \mu(\mathrm{d}y) = \infty$ (by \cite[Theorem 10.4.10]{2009-M-T}) and the equivalence of $\mu$ and $\psi$. $\Box$

\textbf{Proof of necessity of Theorem \ref{Non-stro}} Since the chain $\{X_n\}$ is non-strongly ergodic, there exists a set $A \in \mathscr{B}^+(\mathscr{X})$ such that $\sup_{x \in \mathscr{X}} \mathbb{E}_x \tau_A^+ = \infty$. Without loss of generality, we may assume $A$ is compact (otherwise, by \cite[Theorem 16.0.2, Proposition 6.2.8]{2009-M-T}, $\sup_{x \in \mathscr{X}} \mathbb{E}_x \tau_{E_m}^+ = \infty$ for some $m \ge 1$, so we can replace $A$ by the closure of a suitable set).

Since $\{P^{(m)}(x, \cdot)\}$ on $(E_m, \mathscr{B}(E_m))$ is strongly ergodic, $\sup_{x \in E_m \cap A^c} V_A^{*(m)}(x) < \infty$, where $V_A^{*(m)}(x) = \mathbb{E}_x^{(m)} \tau^{(m)+}_A$. Define $W^{(m)}(x) := \mathbf{1}_{E_m \cap A^c}(x) V_A^{*(m)}(x)$ for all $x \in \mathscr{X}$.
 It is obvious that $\sup_{x \in A^c} W^{(m)}(x) < \infty$ for all $m \ge 1$.
For $x \in E_m \cap A^c$, $P^{(m)}(x, \cdot) = P(x, \cdot)$, so by Proposition \ref{algbra-Moments},
\[
\begin{aligned}
W^{(m)}(x) &= \int_{A^c \cap E_m} V_A^{*(m)}(y) P^{(m)}(x, \mathrm{d}y) + 1 \\
&= \int_{A^c \cap E_m} W^{(m)}(y) P(x, \mathrm{d}y) + 1.
\end{aligned}
\]
This implies $W^{(m)}(x) \le P W^{(m)}(x) + 1$ for all $x \in A^c$.

By Lemma \ref{lem-nece}, $\lim_{m \to \infty} \uparrow W^{(m)}(x) = \mathbb{E}_x \tau_A^+$ for all $x \in A^c$. Hence, $\sup_{\substack{x \in A^c \\ m \ge 1}} W^{(m)}(x) = \infty$. $\Box$

\subsubsection{Proof of Theorem  \ref{Non-geo}}
%The following lemmas are needed for the proof of Theorem \ref{Non-geo}:
\begin{mylem}\label{non-geo-lem1} For a recurrent Markov chain $\{X_n\}$, a set $A \in \mathscr{B}(\mathscr{X})$, and $r > 1$, $\mathbb{E}_x \sum_{n=0}^{\tau_A^+ - 1} r^n$ is the minimal non-negative solution to
\[
V(x) = r \int_{A^c} V(y) P(x, \mathrm{d}y) + 1, \quad x \in \mathscr{X}.
\]
\end{mylem}

\textbf{Proof} By \cite[Corollary 2.8]{2014-M-S}, $\mathbb{E}_x r^{\tau_A^+} \mathbf{1}_{\{\tau_A^+ < \infty\}}$ is the minimal solution to
\begin{equation}\label{F-L-1}
V(x) = r \int_{A^c} V(y) P(x, \mathrm{d}y) + r P(x, A), \quad x \in \mathscr{X}.
\end{equation}
Since $\{X_n\}$ is recurrent, $\mathbb{P}_x(\tau_A^+ < \infty) = 1$ for all $x \in \mathscr{X}$. Using the geometric series identity
\[
\mathbb{E}_x \sum_{n=0}^{\tau_A^+ - 1} r^n = \frac{1}{r - 1} \mathbb{E}_x \big( r^{\tau_A^+} - 1 \big),
\]
the required assertion follows directly from \eqref{F-L-1}. $\Box$

\iffalse

%\begin{mylem}\label{non-geo-lem} For a $\psi$-irreducible and aperiodic Markov chain $\{X_n\}$, if there exist sets $A, B \in \mathscr{B}^+(\mathscr{X})$ such that for all $r > 1$,
\[
\left\{ x \in B : \mathbb{E}_x r^{\tau_A^+} = \infty \right\} \in \mathscr{B}^+(\mathscr{X}),
\]
then $\{X_n\}$ is not geometrically ergodic. 
\end{mylem}
%\textbf{Proof} By \cite[Propositions 15.3.1 and 15.3.2]{2009-M-T}, the assertion holds immediately. $\Box$
\fi

Let $P^{(m)}(x, \cdot)$ be a Markov kernel on $(E_m, \mathscr{B}(E_m))$ defined in Proposition \ref{approx}. Then $P^{(m)}(x, \cdot) = P(x, \cdot)$ for all $x \in E_m \cap A^c$. Consider the following equation for $r > 1$:
\[
V(x) = r \int_{A^c \cap E_m} V(y) P(x, \mathrm{d}y) + 1, \quad x \in E_m.
\]
Denote its minimal solution by $V_A^{*(r, E_m)}(x)$ for $x \in E_m$. By Lemma \ref{non-geo-lem1} and Lemma \ref{Monot}, we have
\begin{equation}\label{geo-eq}
V_A^{*(r, E_m)}(x) \uparrow \mathbb{E}_x \sum_{k=0}^{\tau_A^+ - 1} r^k \quad \text{as } m \to \infty \text{ for all } x \in A.
\end{equation}

\textbf{Proof of Theorem \ref{Non-geo}} We first show that
\[
V^{(n)}(x) \le \mathbb{E}_x \sum_{k=0}^{\tau_A^+ - 1} r_n^k \quad \text{for all } x \in A.
\]

Since $V^{(n)}(x)$ is compactly supported for each $n \ge 1$, we may choose an increasing sequence of compact sets $E_n \uparrow \mathscr{X}$ such that
\[
V^{(n)}(x) \le r_n \int_{A^c \cap E_n} V^{(n)}(y) P(x, \mathrm{d}y) + 1 \quad \text{for all } n \ge 1 \text{ and } x \in E_n \cap A^c.
\]

Let $W^{(r_n, E_n)}(x)$ for $x \in E_n$ denote the minimal solution to
\[
W(x) = r_n \int_{A^c \cap E_n} W(y) P(x, \mathrm{d}y) + 1 = r_n \int_{A^c \cap E_n} W(y) P^{(n)}(x, \mathrm{d}y) + 1
\]
for all $x \in E_n \cap A^c$ and $n \ge 1$. By Lemma \ref{Second-ite-geo} and the fact that $\sup_{x \in E_n} V^{(n)}(x) < \infty$, we obtain
\[
V^{(n)}(x) \le W^{(r_n, E_n)}(x) \quad \text{for all } x \in E_n \cap A^c.
\]

Moreover, by \eqref{geo-eq}, for all sufficiently large $n$,
\[
V^{(n)}(x) \le W^{(r_n, E_n)}(x) \le \mathbb{E}_x \sum_{k=0}^{\tau_A^+ - 1} r_n^k \quad \text{for all } x \in E_n \cap A^c.
\]

For $x \in A$, we thus have
\[
\begin{aligned}
V^{(n)}(x) &\le r_n \int_{A^c \cap E_n} \mathbb{E}_y \sum_{k=0}^{\tau_A^+ - 1} r_n^k P(x, \mathrm{d}y) + 1 \\
&\le r_n \int_{A^c} \mathbb{E}_y \sum_{k=0}^{\tau_A^+ - 1} r_n^k P(x, \mathrm{d}y) + 1 \\
&= \mathbb{E}_x \sum_{k=0}^{\tau_A^+ - 1} r_n^k.
\end{aligned}
\]

For any $r > r_n$, it follows that
\[
V^{(n)}(x) \le \mathbb{E}_x \sum_{k=0}^{\tau_A^+ - 1} r_n^k \le \mathbb{E}_x \sum_{k=0}^{\tau_A^+ - 1} r^k.
\]

By condition (3) of Theorem \ref{Non-geo}, the Markov chain is non-geometrically ergodic. $\Box$

\section{A class of Markov process with immigration}
In this section, we investigate a type of examples to illustrate our criteria.

\subsection{Example 1}
Our criteria  are especially effective for models like the following one.

Let $1\ge b(x)>0$ and $1\ge\gamma(x)>0$ be two continuous functions on $\mathbb{R}^+$ with $b(x)+\gamma(x)\le1$. Let $\beta(x)$ be a density function with support $(0, a+1)$ ($a\le \infty$). $(X_n)_{n\ge0}$ is a Markov chain on $\mathbb{R}^+$ with one-step transition kernel $P(x, {\rm d}y)$ given by
\begin{equation*}
P(x,{\rm d}y)=\begin{cases}
\beta(y){\rm d}y,&x=0\\
\gamma(x)\delta_{\{0\}}({\rm d}y)+b(x)\delta_{\{x+1\}}({\rm d}y)+(1-\gamma(x)-b(x))\delta_{\{x\}}({\rm d}y), &x>0.
\end{cases}
\end{equation*}
\iffalse
\begin{center}
	\begin{tikzpicture}[>=latex]
		% 绘制数轴
	\draw[->, thick] (0,0) -- (10,0);
	\foreach \x in {0,4,5,9} {
		\draw (\x,0.1) -- (\x,-0.1); % 绘制刻度线
	}
			
	% 标记点
	\node at (0,-0.3) {0};
	%\node at (1,-0.3) {1};
	%\node at (3,-0.3) {$x-1$};
	\node at (4,-0.3) {$x$};
	\node at (5,-0.3) {$x+1$};
	%\node at (6,-0.3) {$x+2$};
	\node at (9,-0.3) {$y$};
	\draw[red, thick, ->] (4,0) to[out=-150, in=-30, looseness=1]  node[below, pos=0.6, sloped] {$\gamma(x)$}(0,0);
	%\draw[black, thick, ->] (4,0) to[out=120, in=60, looseness=2] (3,0);
	\draw[black, thick, ->] (4,0) to[out=60, in=120, looseness=2]  node[above, pos=0.2, sloped] {$b(x)$}(5,0);
\draw[blue, thick, ->] (0,0) to[out=15, in=165, looseness=1.5] node[above, pos=0.3, sloped] {$\beta(y)$}(9,0);
%\draw[blue, thick, ->] (0,0) to[out=15, in=165, looseness=1.5] (1,0);
\draw[blue, thick, ->] (0,0) to[out=15, in=165, looseness=1.5]node[below, pos=0.5, sloped] {$\beta(x)$}(4,0);
	\end{tikzpicture}
	\end{center}
	\fi
It is obvious that $\{X_n\}$ is aperiodic and irreducible with respect to Lebesgue measure $m({\rm d}x)$. The positivity of the functions is intended to simplify the discussion, although we can weaken the requirement to the functions being non-negative everywhere while preserving the Markov chain's aperiodicity and irreducibility.

By \cite[Theorem 6.2.8]{2009-M-T}, all compact subsets of $\mathcal{X}$ are petite.

\begin{myprop}\label{ex1-transi} If there exists $x_1>0$ such that $b(x)>0$ for $x\ge x_1$, and
$$\prod_{m=[x_1]}^{\infty}\inf_{m\le x<m+1}\frac{b(x)+\gamma(x)}{b(x)}=\infty,$$
then $\{X_n\}$ is recurrent. Conversely, if there exists $x_1>0$ such that $b(x)>0$ for $x\ge x_1$, and
$$\prod_{m=[x_1]}^{\infty}\sup_{m\le x<m+1}\frac{b(x)+\gamma(x)}{b(x)}<\infty,$$
then $\{X_n\}$ is transient.
\end{myprop}

\textbf{Proof } (1) Let $C=[0,x_1]$ for some $x_1>0$. Define
\begin{equation*}
V(x)=\inf_{ [x_1]\le y< [x_1]+1}\frac{b(y)+\gamma(y)}{b(y)}, \quad 0\le x < x_1+1; 
\end{equation*}
\begin{equation*}
V(x)=\prod_{k=[x_1]}^{[x]-1}\inf_{k\le y< k+1}\frac{b(y)+\gamma(y)}{b(y)}, \quad x\ge x_1+1.
\end{equation*}
Then 
$$V(x+1)= \inf_{[x]\le y<[x]+1}\frac{b(y)+\gamma(y)}{b(y)}V(x)\le \frac{b(x)+\gamma(x)}{b(x)}V(x),\quad x>x_1,$$
i.e., $PV(x)-V(x)\le0$ for $x>x_1$. Since
$$\aligned C_V(n)=&\{y: V(y)\le n\}\\=&\left\{y: y\ge x_1+1, \prod_{k=[x_1]}^{[y]-1}\inf_{k\le x<k+1}\frac{b(x)+\gamma(x)}{b(x)}\le n\right\}\\
&\cap\left\{y: 0<y\le x_1+1,  \inf_{[x_1]\le x<[x_1]+1}\frac{b(x)+\gamma(x)}{b(x)}\le n\right\}, \endaligned$$
which is a compact set when $\gamma(x)$ and $b(x)$ are continuous and 
$$\prod_{k=[x_1]}^{\infty}\inf_{k\le x<k+1}\frac{b(x)+\gamma(x)}{b(x)}=\infty.$$ Hence $\{X_n\}$ is recurrent by \cite[Theorem 8.4.3]{2009-M-T}. 

(2) Let $C=[0, x_1]$. Define
\begin{equation*}
V(x)=\begin{cases}
0,&0\le x<x_1\\
-1,&x_1\le x<x_1+1\\
-\sup_{[x_1]\le y<[x_1]+1}\frac{b(y)+\gamma(y)}{b(y)}, & x_1+1\le x< x_1+2,\\
-\prod_{k=[x_1]}^{[x]-1}\sup_{k\le y<k+1}\frac{b(y)+\gamma(y)}{b(y)},&  x\ge x_1+2.
\end{cases}
\end{equation*}
Then $$V(x+1)=-\sup_{[x]\le y< [x]+1}\frac{b(y)+\gamma(y)}{b(y)}V(x)\le-\frac{b(x)+\gamma(x)}{b(x)}V(x),\quad x>x_1,$$
i.e., $PV(x)-V(x)\le0$ for $x>x_1$. Since
$\{x: V(x)<\inf_{y\in C} V(y)\}=C^c\in\mathcal{B}^+(\mathcal{X})$ and 
$$\inf_{y\in\mathcal{X}} \left(-\prod_{k=[x_1]}^{[y]-1}\sup_{k\le x<k+1}\frac{b(x)+\gamma(x)}{b(x)}\right)>-\infty,$$
when $\gamma(x)$ and $b(x)$ are continuous and 
$$\prod_{k=[x_1]}^{\infty}\sup_{k\le x<k+1}\frac{b(x)+\gamma(x)}{b(x)}<\infty.$$ Hence $\{X_n\}$ is transient by Proposition \ref{transient-1}.

\begin{myprop}\label{ex1-erg}
(1) $(X_n)$ is non-ergodic there exists $x_1>0$ such that $$\int_{x_1}^{\infty}\frac{\beta(y)}{\sup_{u\ge y}\gamma(u)}{\rm d}y=\infty.$$

(2) $(X_n)$ is ergodic if there exists $x_1>0$ and $c>0$ such that $$\int_{x_1}^{\infty}\frac{cb(y)+1}{\gamma(y)}\beta(y){\rm d}y<\infty$$ and  
$$b(x)\left(c\left(\frac{b(x+1)}{\gamma(x+1)}-\frac{b(x)}{\gamma(x)}\right)- \frac{1}{\gamma(x)}+\frac{1}{\gamma(x+1)}-c\right) \le0,\quad x>x_1.$$
\end{myprop}
{\bf Proof } (1) Since $(X_n)$ is irreducible with respect to Lebesgue measure, 
we have
$$0<\sup_{y\ge x}\gamma(y)<\infty,\quad x\ge 0.$$
 Let  \begin{equation*}
V^{(n)}(x)=\frac{1}{\sup_{y\ge x\wedge n}\gamma(y)},\quad x>x_1,\qquad
V^{(n)}(x)=\int_{x_1}^{n}\frac{\beta(u)}{\sup_{t\ge u\wedge n}\gamma(t)}{\rm d}t,\quad 0\le x\le x_1.
\end{equation*} Then  
$$\aligned & \int_{x_1}^{\infty}V^{(n)}(y)P(0, {\rm d}y)-V^{(n)}(0)+1=\int_{n}^{\infty}\frac{\beta(y)}{\sup_{t\ge y\wedge n}\gamma(t)}{\rm d}y+1\ge0;\\
&\int_{x_1}^{\infty}V^{(n)}(y)P(x, {\rm d}y)-V^{(n)}(x)+1=b(x)\left(V(x+1)-V(x)\right)-\gamma(x)V(x)+1\ge0,\quad x>0.\endaligned$$
%$$\int_0^{x_1}\sup_{n\ge1}V^{(n)}(y){\rm d}y=x_1\int_{x_1}^{\infty}\frac{\beta(u)}{\sup_{t\ge u}\gamma(t)}{\rm d}u=\infty.$$
Hence,  $(X_n)$ is non-ergodic by letting $A=[0, x_1]$ in Corollary \ref{Non-1-ergo-cor}. %Theorem \ref{Non-1-ergo}.

(2) Let $V(x)=0$ for $0\le x\le x_1$, \begin{equation*}
V(x)=\frac{cb(x)+1}{\gamma(x)},\quad x>x_1.
\end{equation*} Then $PV(x)-V(x)+1\le0$ for $x>x_1$
is equivalent to 
$$b(x)\left(c\left(\frac{b(x+1)}{\gamma(x+1)}-\frac{b(x)}{\gamma(x)}\right)- \frac{1}{\gamma(x)}+\frac{1}{\gamma(x+1)}-c\right)\le0,\quad x>x_1.$$
The required assertion holds by \cite[Theorem 13.0.1]{2009-M-T}.

\begin{myprop}\label{ex1-str-erg}
(1) $(X_n)$ is non-strongly ergodic if $\overline{\lim}_{x\to\infty}\gamma(x)=0$.

(2) If $\underline{\lim}_{x\to\infty}\gamma(x)>0$, then $(X_n)$ is strongly ergodic.  

(3) If $\overline{\lim}_{x\to\infty}\gamma(x)=0$, $\underline{\lim}_{x\to\infty}b(x)>0$, then $\{X_n\}$ is non-geometrically ergodic.
\end{myprop}

{\bf Proof} (1)
Since $\{X_n\}$ is irreducible with respect to Lebesgue measure, we get
$0<\sup_{y\ge x}\gamma(y)<\infty$ for $x\ge0$.
Let $V^{(n)}(x)=\left(\sup_{y\ge x\wedge n}\gamma(y)\right)^{-1}{\bf 1}_{[0,x_1]}(x)$ for some $x_1>0$.
Then $(X_n)$ is non-strongly ergodic if
$\overline{\lim}_{x\to\infty}\gamma(x)=0$ by Theorem \ref{Non-stro}.

(2) By $\underline{\lim}_{x\to\infty}\gamma(x)>0$, there exists some $x_1>0$ such that
$$0<\inf_{y\ge x}\gamma(y),\quad x\ge x_1.$$
Let $\beta:=\inf_{y\ge x_1}\gamma(y)>0$, \begin{equation}\label{zsy-1}
V(x)=\left({\inf_{y\ge x\vee (x_1+1)}\gamma(y)}\right)^{-1}.
\end{equation} 
Then $V(x)\ge1$, $0<\beta\le1$ and
$$\aligned 
&\int_{0}^{\infty}V(y)P(x, {\rm d}y)+(\beta-1)V(x)\le 0,\quad x>0,\endaligned$$
$$\aligned 
\int_{0}^{\infty}V(y)P(0, {\rm d}y)+(\beta-1)V(0)&\le \frac{1}{\inf_{u\ge x_1+1}\gamma(u)}\int_{0}^{\infty}\beta(y){\rm d}y<\infty.\endaligned$$
By \cite[Theorem 16.0.2]{2009-M-T}, $(X_n)$ is strongly ergodic if $\underline{\lim}_{x\to\infty}\gamma(x)>0$.

\iffalse
(4) Let $V(x)$ and $\beta(x)$ be the same as those in (3) above. 
By \cite[Theorem 15.0.1]{2009-M-T}, $\{X_n\}$ is geometrically ergodic if there exists $x_1>0$ such that $$\int_{x_1+1}^{\infty}\beta(y)/\inf_{u\ge y}\gamma(u){\rm d}y<\infty$$ and $\inf_{x\ge y_0}\gamma(x)>0$ for some $y_0$.
\fi

(3) If $\{X_n\}$ is geometrically ergodic, by Lemma \ref{non-geo-lem1}, the following equation has a finite nonnegative solution for some $A\in\mathcal{B}^+(\mathcal{X})$ and $r>1$.
\begin{equation}\label{non-geo-lem-ap1}
V(x) = r \int_{A^c} V(y) P(x, \mathrm{d}y) + 1, \quad x \in \mathscr{X}.
\end{equation}
Without loss of generality, assume that $A=[0,1]$. By
\eqref{non-geo-lem-ap1},
\begin{equation}\label{non-geo-lem-ap2}V(x)=r\left[b(x)V(x+1)+(1-b(x)-\gamma(x))V(x)\right]+1,\quad x>0.\end{equation}
Moreover, for $x>x_1$, $x_1>0$,
\begin{equation*}
\begin{split}
V(x+1) - V(x) &= \frac{\gamma(x) + 1/r - 1}{b(x)}V(x) - \frac{1/r}{b(x)} \\
&\le \frac{V(x)}{\inf_{x > x_1} b(x)} \left( \sup_{x \ge x_1} \gamma(x) + \frac{1}{r} - 1 \right).
\end{split}
\end{equation*}
Since $\overline{\lim}_{x\to\infty}\gamma(x)=0$, $\underline{\lim}_{x\to\infty}b(x)>0$, we have $V(x+1)\le V(x)$ for large enough $x$ (say $x_2>x_1$).
So $\{V(x): x\ge x_2\}$ is bounded. Since $b(x)$, $\gamma(x)$ are continuous, by \eqref{non-geo-lem-ap2},
$$V(x)=\frac{\gamma(x)+b(x)+1/r}{b(x)}V(x+1)-\frac{1/r}{b(x)},$$
we can get $\{V(x): x\in A^c\cap[0, x_2]\}$ is bounded. Furthermore, $\{V(x): x\ge0\}$ is bounded by \eqref{non-geo-lem-ap1}.
Consequently, $\left(\mathbb{E}_x\sum_{n=0}^{\tau_A^+-1}r^n: x\ge 0\right)$ is bounded since it is the minimal non-negative solution to \eqref{non-geo-lem-ap1}, so is $(\mathbb{E}_x\tau^+_A: x\ge 0)$. 
$\{X_n\}$ is therefore strongly ergodic. But this is impossible. Hence, $\{X_n\}$ is non-geometrically ergodic. 

\textbf{Proof of Proposition \ref{1-Ex-1}} This is a special case of $b(x)+\gamma(x)=1$ in Proposition \ref{ex1-transi} $\sim$ Proposition \ref{ex1-str-erg}. We only mention that when proving the ergodicity of the Markov chain, since $1>\gamma(x)>0$,
\[
\int_{x_1}^{\infty}\frac{c(1-\gamma(x))+1}{\gamma(y)}\beta(y)\,\mathrm{d}y\le
(c+1)\int_{x_1}^{\infty}\frac{\beta(y)}{\gamma(y)}\,\mathrm{d}y,
\]
and
\[
(1-\gamma(x))\left(c\left(\frac{1-\gamma(x+1)}{\gamma(x+1)}-\frac{1-\gamma(x)}{\gamma(x)}\right)- \frac{1}{\gamma(x)}+\frac{1}{\gamma(x+1)}-c\right) \le0,\quad x>x_1
\]
is equivalent to
\[
\frac{1}{\gamma(x+1)}-\frac{1}{\gamma(x)}\le \frac{c}{c+1}<1,\quad x>x_1.
\]

\subsection{Example 2}
We investigate the instability of a special type of single death process with immigration defined on $\mathbb{Z}^+$.

Suppose that the nonnegative sequences $\{\beta_i\}$, $\{\gamma_i\}$, and $\{p_i\}$ satisfy $0 < \gamma_i + p_i \le 1$ for $i = 2, 3, \cdots$, and $\sum_{j=0}^{\infty}\beta_j = 1$. The Markov chain $(X_n)_{n\ge0}$ is irreducible with respect to the counting measure, and its one-step transition probability $P_{ij}$ is given by:
\begin{equation}\label{trans-ex}
 P_{ij}=\begin{cases}
\beta_j,&i=0,\,j\ge0,\\
\gamma_i, &i\ge2,\, j=0,\\
p_i, & i\ge1,\, j=i-1,\\
1-p_i-\gamma_i{\bf 1}_{\{i\ge 2\}},& j=i\ge1.
\end{cases}
\end{equation}

\begin{center}
	\begin{tikzpicture}[>=latex]
		% Draw the number line
		\draw[->, thick] (0,0) -- (10,0);
		\foreach \x in {0,1,3,4,5,6,9} {
			\draw (\x,0.1) -- (\x,-0.1); % Draw tick marks
		}
			
		% Label points
		\node at (0,-0.3) {0};
		\node at (1,-0.3) {1};
		\node at (3,-0.3) {$i-1$};
		\node at (4,-0.3) {$i$};
		\node at (5,-0.3) {$i+1$};
		\node at (6,-0.3) {$i+2$};
		\node at (9,-0.3) {$i+k$};
		\draw[red, thick, ->] (4,0) to[out=-150, in=-30, looseness=1] (0,0);
		\draw[red, thick, ->] (5,0) to[out=-150, in=-30, looseness=1] (0,0);
		\draw[black, thick, ->] (5,0) to[out=120, in=60, looseness=2] (4,0);
		\draw[black, thick, ->] (6,0) to[out=120, in=60, looseness=2] (5,0);
		\draw[black, thick, ->] (4,0) to[out=120, in=60, looseness=2] (3,0);
		\draw[blue, thick, ->] (0,0) to[out=15, in=165, looseness=1.5] (9,0);
		\draw[blue, thick, ->] (0,0) to[out=15, in=165, looseness=1.5] (1,0);
		\draw[blue, thick, ->] (0,0) to[out=15, in=165, looseness=1.5] (3,0);
	\end{tikzpicture}
\end{center}

Since $(X_n)$ is irreducible, it must hold that $\gamma_i + p_i > 0$ for all $i \ge 2$ and $p_1 > 0$. The Markov chain is clearly recurrent at all times. In fact, let $v_i = i$ for $i \ge 0$; this function is finite and satisfies $\lim_{i\to\infty}v_i = \infty$, and for any $i \ge 1$, we have
\[
\sum_{j=0}^{\infty}P_{ij}v_j \le v_i.
\]
By Theorem 4.24 in \cite{2004-Chen}, the above recurrence conclusion is proven.

For the criteria of instability, we have the following results:

\begin{myprop}\label{Prop2} Let the Markov chain $\{X_n\}$ be generated by the one-step transition matrix \eqref{trans-ex}. Then:
\begin{itemize}
\item[(1)] $(X_n)$ is non-ergodic if there exists $i_0 > 0$ such that
\[
\sum_{j=i_0}^{\infty}\frac{\beta_j}{p_j + \gamma_j} = \infty;
\]
\item[(2)] $(X_n)$ is non-geometrically ergodic if there exists $a > 0$ such that $\gamma_i + ap_i \le i^{-a}$ for all $i \ge 2$;
\item[(3)] $(X_n)$ is non-strongly ergodic if $\inf_{i\ge2}(p_i + \gamma_i) = 0$.
\end{itemize}
\end{myprop}

\noindent{\bf Proof} 
(1) Let $A = \{0\}$ and define the auxiliary function as
\[
v_1^{(n)} = \frac{1}{p_1}, \quad v_i^{(n)} = \frac{1}{p_i + \gamma_i}\quad (i \ge 2),\quad v_0^{(n)} = \sum_{k=2}^{n}\frac{\beta_k}{p_k + \gamma_k}.
\]
By Theorem \ref{Non-1-ergo}, $(X_n)$ is non-ergodic if $\sum_{j=i_0}^{\infty}\frac{\beta_j}{p_j + \gamma_j} = \infty$.

(3) Let $v_i^{(n)} = i^a{\bf 1}_{\{i\le n\}}$ ($a > 0$) and $r_n = 1 + \frac{1}{n}$. It is easy to verify that for $i = 0, 1$ or $i \ge n$,
\[
r_n\sum_{j=1}^{\infty}P_{ij}v_j^{(n)} - v_i^{(n)} + 1  \ge 0.
\]
For $n > i \ge 2$, using $i^a - (i-1)^a \le ai^a$ and $r_n \ge 1$, we obtain
\[
\begin{aligned}
r_n\sum_{j=1}^{\infty}P_{ij}v_j^{(n)} - v_i^{(n)} + 1 &= r_n (1 - p_{i} - \gamma_i)i^a + r_np_{i}(i-1)^a - i^a + 1 \\
&\ge -p_iai^a - \gamma_i i^a + 1 \ge 0.
\end{aligned}
\]
By Theorem \ref{Non-geo}, $(X_n)$ is non-geometrically ergodic.

(4) Let $A = \{0, 1\}$ and define the auxiliary function as $v_i^{(n)} = 0$ for $i = 0, 1$, and
\[
v_i^{(n)} = \frac{1}{p_i + \gamma_i}\quad (i \ge 2).
\]
Then
\[
\begin{aligned}
Pv^{(n)}(i) - v_i^{(n)} + 1 &= -(p_i + \gamma_i)v_{i}^{(n)} + p_iv_{i-1}^{(n)} + 1 \\
&\ge 0\quad (i \ge 2).
\end{aligned}
\]
By Theorem \ref{Non-stro}, $(X_n)$ is non-strongly ergodic if $\inf_{j\ge 2}(p_j + \gamma_j) = 0$.

\subsection*{Acknowledgments}
This work is supported by the National Key R\&D Program of China (No. 2020YFA0712900) and the National Natural Science Foundation of China (No. 12101186, No. 12171038).


\begin{thebibliography}{99}
\setlength{\itemsep}{-2mm}

\bibitem{2004-C-K}
B.D. Choi, B. Kim, Non-ergodicity criteria for denumerable continuous-time Markov processes, \emph{Operations Research Letters} 32(6) (2004) 574--580.

\bibitem{2004-Chen}
M.-F. Chen, \emph{From Markov Chains to Non-Equilibrium Particle Systems}, 2nd edition, World Scientific, Singapore, 2004.

\bibitem{2002-J-R}
S.F. Jarner, G.O. Roberts, Polynomial convergence rates of Markov chains, \emph{Annals of Applied Probability} 12(1) (2002) 224--247.

\bibitem{1995-L}
O.S. Lee, Geometric ergodicity and transience for nonlinear autoregressive models, \emph{Communications of the Korean Mathematical Society} 10(2) (1995) 409--417.

\bibitem{2003-Mao}
Y.-H. Mao, Algebraic convergence for discrete-time ergodic Markov chains, \emph{Sci. China Ser. A} 46(5) (2003) 621--630.

\bibitem{2021-M-W}
Y.-H. Mao, T. Wang, Lyapunov-type conditions for non-strong ergodicity of Markov processes, \emph{Journal of Applied Probability} 58 (2021) 238--253.

\bibitem{2014-M-S}
Y.-H. Mao, Y.-H. Song, On geometric and algebraic transience for discrete-time Markov chains, \emph{Stochastic Processes and Their Applications} 124(4) (2014) 1648--1678.

\bibitem{1992-M-T}
S.P. Meyn, R.L. Tweedie, Stability of Markovian processes I: criteria for discrete-time chains, \emph{Advances in Applied Probability} 24(3) (1992) 542--574.

\bibitem{1993-M-T}
S.P. Meyn, R.L. Tweedie, Stability of Markovian Processes III: Foster-Lyapunov criteria for continuous-time processes, \emph{Advances in Applied Probability} 25(3) (1993) 518--548.

\bibitem{2009-M-T}
S.P. Meyn, R.L. Tweedie, \emph{Markov Chains and Stochastic Stability}, 2nd edition, Cambridge University Press, 2009.

\bibitem{1983-S-T}
L.I. Sennott, H.R.L. Tweedie, Mean drifts and the non-ergodicity of Markov chains, \emph{Operations Research} 31(4) (1983) 783--789.

\bibitem{1985-S}
W. Szpankowski, Some sufficient conditions for non-ergodicity of Markov chains, \emph{Journal of Applied Probability} 22(1) (1985) 138--147.

\bibitem{1990-T}
D. Tj{\o}stheim, Non-linear time series and Markov chains, \emph{Advances in Applied Probability} 22 (1990) 587--611.

\bibitem{2022-Wei}
Z.F. Wei, Inverse problems for ergodicity of Markov chains, \emph{Journal of Mathematical Analysis and Applications} 505 (2022) 125483.

\bibitem{2010-Zhang-Zhao}
Y.H. Zhang, Q.Q. Zhao, \emph{Single Birth Processes with Immigration}, Acta Mathematica Sinica (English Series), Vol. 26, No. 5, pp. 833--846, 2010.

\end{thebibliography}
\end{document}